\documentclass{nyj}

\newcommand\R{{\hbox{\bf R}}}

\renewcommand\Im{{\hbox{Im}}}
\renewcommand\Re{{\hbox{Re}}}
\newcommand\Mass{{\hbox{\rm Mass}}}
\newcommand\dist{{\hbox{\rm dist}}}

\newcommand\C{{\hbox{\bf C}}}

\newcommand\eps{\varepsilon}

\theoremstyle{plain}
  \newtheorem{theorem}[subsection]{Theorem}
  
  \newtheorem{proposition}[subsection]{Proposition}
  \newtheorem{lemma}[subsection]{Lemma}
  \newtheorem{corollary}[subsection]{Corollary}

\theoremstyle{remark}
  \newtheorem{remark}[subsection]{Remark}

\theoremstyle{definition}

\include{psfig}

\begin{document}

\title[Scattering for energy-critical NLS]{Global well-posedness and scattering for the higher-dimensional
energy-critical non-linear Schr\"odinger equation for radial data}
\author{Terence Tao}
\address{Department of Mathematics, UCLA, Los Angeles CA 90095-1555}
\email{tao@@math.ucla.edu}
\thanks{The author is a Clay Prize Fellow and is supported by the Packard
Foundation.  The author is indebted to Jean Bourgain, Jim Colliander, Manoussos Grillakis,
Markus Keel, Gigliola Staffilani, and Hideo Takaoka for useful conversations.  The author
also thanks Monica Visan and the anonymous referee for several corrections.}

\vspace{-0.3in}
\begin{abstract}
In any dimension $n \geq 3$, we show that spherically symmetric bounded energy solutions 
of the defocusing energy-critical non-linear Schr\"odinger 
equation $i u_t + \Delta u = |u|^{\frac{4}{n-2}} u$ in $\R \times \R^n$ exist globally and scatter to free solutions;
this generalizes the three and four dimensional results of Bourgain \cite{borg:scatter}, \cite{borg:book}
and Grillakis \cite{grillakis:scatter}.
Furthermore we have bounds on various spacetime norms of the solution which are of exponential type in the energy,
which improves on the tower-type bounds of Bourgain.  In higher dimensions $n \geq 6$ some new technical difficulties arise because of the very low power of the non-linearity.
\end{abstract}

\maketitle

\section{Introduction}

Let $n \geq 3$ be an integer.  We consider solutions $u: I \times \R^n \to \C$ of the
defocusing energy-critial non-linear Schr\"odinger equation
\begin{equation}\label{nls}
i u_t + \Delta u = F(u)
\end{equation}
on a (possibly infinite) time interval $I$, where $F(u) := |u|^{\frac{4}{n-2}} u$.  
We will be interested in the Cauchy problem
for the equation \eqref{nls}, specifying initial data $u(t_0)$ for some $t_0 \in I$ and then studying the existence 
and long-time behavior of solutions to this Cauchy problem.  

We restrict our attention to solutions for which the energy
$$ E(u) = E(u(t)) := \int_{\R^n} \frac{1}{2} |\nabla u(t,x)|^2 + \frac{n-2}{2n} |u(t,x)|^{\frac{2n}{n-2}}\ dx$$
is finite.  It is then known (see e.g. \cite{cwI}) that for any given choice of finite energy initial data $u(t_0)$,
the solution exists for times close to $t_0$, and the energy $E(u)$ is conserved in those times.  Furthermore
this solution is unique\footnote{In fact, the condition that the solution lie in $L^{2(n+2)/(n-2)}_{t,x}$ can
be omitted from the uniqueness result, thanks to the endpoint Strichartz estimate in \cite{tao:keel} and
the Sobolev embedding $\dot H^1_x \subseteq L^{2n/(n-2)}_x$; see \cite{kato}, \cite{ft}, \cite{ftp} for further discussion.
We thank Thierry Cazenave for this observation.}
 in the class $C^0_t \dot H^1_x \cap L^{2(n+2)/(n-2)}_{t,x}$, and we shall always assume
our solutions to lie in this class.
The significance of the exponent
in \eqref{nls} is that it is the unique exponent which is \emph{energy-critical}, in the sense that
the natural scale invariance
\begin{equation}\label{scaling}
 u(t,x) \mapsto \lambda^{-(n-2)/2} u(\frac{t}{\lambda^2}, \frac{x}{\lambda})
\end{equation}
of the equation \eqref{nls} leaves the energy invariant; in other words, the energy $E(u)$ is a dimensionless
quantity.

If the energy $E(u(t_0))$ is sufficiently small (smaller than some absolute constant 
$\eps > 0$ depending only on $n$) then it is known (see \cite{cwI}) that one has a unique global finite-energy
solution $u: \R \times \R^n \to \C$ to \eqref{nls}.  Furthermore we have the global-in-time 
\emph{Strichartz bounds}
$$ \| \nabla u \|_{L^q_t L^r_x(\R \times \R^n)} \leq C(q,r,n,E(u))$$
for all exponents $(q,r)$ which are \emph{admissible} in the sense that\footnote{Strictly speaking, the
result in \cite{cwI} did not obtain these estimates for the endpoint $q=2$, but they can easily be recovered by
inserting the Strichartz estimates from \cite{tao:keel} into the argument in \cite{cwI}.}
\begin{equation}\label{admissible}
2 \leq q,r \leq \infty; \quad \frac{1}{q} + \frac{n}{2r} = \frac{n}{4}.
\end{equation}
In particular, from Sobolev embedding we have the spacetime estimate
\begin{equation}\label{u-spacetime}
 \| u \|_{L^{2(n+2)/(n-2)}_{t,x}(\R \times \R^n)} \leq M(n,E(u))
\end{equation}
for some explicit function $M(n,E) > 0$.
Because of this and some further Strichartz analysis, one can also show scattering, in the sense that there exist
Schwarz solutions $u_+, u_-$ to the \emph{free} Schr\"odinger equation $(i\partial_t + \Delta) u_\pm = 0$,
such that
$$ \| u(t) - u_\pm(t) \|_{\dot H^1(\R^n)} \to 0 \hbox{ as } t \to \pm \infty.$$
This can then be used to develop a small energy scattering theory (existence of wave operators, asymptotic
completeness, etc.); see \cite{caz}.  Also, one can show that the solution map $u(t_0) \to u(t)$ extends to a 
globally Lipschitz map in the energy space $\dot H^1(\R^n)$.

The question then arises as to what happens for large energy data.  In \cite{cwI} it was shown that the Cauchy problem is
locally well posed for this class of data, so that we can construct solutions for short times at least; the issue is
whether these solutions can be extended to all times, and whether one can obtain scattering results like before.
It is well known that such results will indeed hold if one could obtain the \emph{a priori} bound \eqref{u-spacetime} for all
global Schwarz solutions $u$ (see e.g. \cite{borg:book}).  It is here that the sign of the non-linearity in \eqref{nls}
is decisive (in contrast to the small energy theory, in which it plays no role).  Indeed, if we replaced the non-linearity $F(u)$
by the focusing non-linearity  $-F(u)$ then an argument of Glassey \cite{glassey} shows that large 
energy Schwarz initial data can blow up in finite time; for instance, this will occur whenever the potential energy
exceeds the kinetic energy.  

In the defocusing case, however, the existence of \emph{Morawetz inequalities} allows one to obtain better control
on the solution.  A typical such inequality is
$$
\int_I \int_{\R^n} \frac{|u(t,x)|^{2n/(n-2)}}{|x|}\ dx dt \leq C (\sup_{t 
\in I} \| u(t) \|_{\dot H^{1/2}(\R^n)})^2$$
for all time intervals $I$ and all Schwarz solutions $u: I \times \R^n \to \C$ to \eqref{nls}, where $C > 0$ is a constant depending only on $n$; 
this inequality can be 
proven by differentiating the quantity $\int_{\R^n} \Im(\frac{x}{|x|} \cdot \nabla u(t,x) \overline{u(t,x)})\ dx$
in time and integrating by parts.  This inequality is not directly useful for the energy-critical problem, as the right-hand side
involves the Sobolev norm $\dot H^{1/2}(\R^n)$ instead of the energy norm $\dot H^1(\R^n)$.  However, by applying
an appropriate spatial cutoff, Bourgain \cite{borg:scatter}, \cite{borg:book} and Grillakis \cite{grillakis:scatter}
obtained the variant Morawetz estimate
\begin{equation}\label{morawetz}
\int_I \int_{|x| \leq A |I|^{1/2}} \frac{|u(t,x)|^{2n/(n-2)}}{|x|}\ dx dt \leq C A |I|^{1/2} E(u)
\end{equation}
for all $A \geq 1$, where $|I|$ denotes the length of the time interval $I$;
this estimate is more useful as it involves the energy on the right-hand side.
For sake of self-containedness we present a proof of this inequality in Section \ref{review-sec}.

The estimate \eqref{morawetz} is useful for preventing concentration of $u(t,x)$ at the spatial origin $x=0$.
This is especially helpful in the \emph{spherically symmetric case} $u(t,x) = u(t,|x|)$, since the spherical
symmetry, combined with the bounded energy assumption can be used to show that $u$ cannot concentrate at any other 
location than the spatial origin.  Note that spatial concentration is the primary obstruction to establishing global existence for
the critical NLS \eqref{nls}; see e.g. \cite{keraani} for some dicussion of this issue.

With the aid of \eqref{morawetz} and several additional arguments, Bourgain \cite{borg:scatter}, 
\cite{borg:book} and Grillakis \cite{grillakis:scatter} were able to show global existence of large energy spherically
smooth solutions in the three dimensional case $n=3$.  Furthermore, the argument in \cite{borg:scatter}, \cite{borg:book}
extends (with some technical difficulties) to the case $n=4$ and also gives the spacetime bound \eqref{u-spacetime}
(which in turn yields the scattering and global well-posedness results mentioned earlier).  However, the dependence of the constant
$M(n,E(u))$ in \eqref{u-spacetime} on the energy $E(u)$ given by this argument is rather poor; in fact it is an iterated
tower of exponentials of height $O(E(u)^C)$.  This is because the argument is based on an \emph{induction on energy}
strategy; for instance when $n=3$ one selects a small number $\eta > 0$ which depends polynomially on the energy, removes a small component
from the solution $u$ to reduce the energy from $E(u)$ to $E(u) - \eta^4$, applies an induction hypothesis asserting
a bound \eqref{u-spacetime} for that reduced solution, and then glues the removed component back in using perturbation 
theory.  The final argument gives a recursive estimate for $M(3, E)$ of the form
$$ M(3, E) \leq C \exp( \eta^C M(3, E - \eta^4)^C )$$
for various absolute constants $C > 0$, and with $\eta = c E^{-C}$.  It is this recursive inequality which yields
the tower growth in $M(3,E)$.  The argument of Grillakis \cite{grillakis:scatter} is not based on an induction on energy,
but is based on obtaining $L^\infty_{t,x}$ control on $u$ rather than Strichartz control (as in \eqref{u-spacetime}), and
it is not clear whether it can be adapted to give a bound on $M(3,E)$.

The main result of this paper is to generalize the result\footnote{We do not obtain regularity results, except
in dimensions $n=3,4$, simply because the non-linearity $|u|^{4/(n-2)} u$ is not smooth in dimensions $n \geq 5$.
Because of this non-smoothness, we will not rely on Fourier-based techniques such as Littlewood-Paley theory,
$X^{s,b}$ spaces, or para-differential calculus, relying instead on the (ordinary) chain rule and some use of H\"older
type estimates.} of Bourgain to general dimensions, and to remove the tower
dependence on $M(n,E)$, although we are still restricted to spherically symmetric data.
As with the argument of Bourgain, a large portion of our argument generalizes to the non-spherically symmetric case; the spherical symmetry is needed only to ensure that the solution concentrates at the spatial origin, and not at 
any other point in spacetime, in order to exploit the Morawetz estimate \eqref{morawetz}.  In light of the recent result
in \cite{gopher} extending the three-dimensional results to general data, it seems in fact likely that at least some of the ideas here can be used in the non-spherically-symmetric setting; see Remark \ref{remark-gopher}.

\begin{theorem}\label{main}  Let $[t_-, t_+]$ be a compact interval, and let $u \in C^0_t \dot H^1([t_-,t_+] \times \R^n)
\cap L^{2(n+2)/(n-2)}_{t,x}([t_-,t_+] \times \R^n)$ be a 
spherically symmetric solution to 
\eqref{nls} with energy $E(u) \leq E$ for some $E > 0$.  Then we have
$$ \| u \|_{L^{2(n+2)/(n-2)}_{t,x}([t_-,t_+] \times \R^n)} \leq C \exp( C E^C )$$
for some absolute constants $C$ depending only on $n$ (and thus independent of $E$, $t_\pm$, $u$).
\end{theorem}

Because the bounds are independent of the length of the time interval $[t_-, t_+]$, it is a standard matter to use this theorem, combined with the local well-posedness theory in \cite{cwI},
to obtain global well-posedness and scattering conclusions for large energy spherically symmetric data; see
\cite{caz}, \cite{borg:book} for details.

Our argument mostly follows that of Bourgain \cite{borg:scatter}, \cite{borg:book}, but avoids the use of
induction on energy using some ideas from other work \cite{grillakis:scatter}, \cite{gopher}, \cite{tao:focusing}.
We sketch the ideas informally as follows.  Following Bourgain, we choose a small parameter $\eta > 0$ depending polynomially on the energy, and then divide the time interval $[t_-, t_+]$ into a finite number of intervals $I_1, \ldots, I_J$, where on each interval the $L^{2(n+2)/(n-2)}_{t,x}$ norm is comparable 
to $c(\eta)$; the task is then to bound the number $J$ of such intervals by $O(\exp(C E^C))$.

An argument of Bourgain based on Strichartz inequalities and harmonic analysis, which we reproduce here\footnote{For some results in the same spirit, showing that ``bubbles'' are the only obstruction to global existence, see \cite{keraani}.},
shows that for each such interval $I_j$, there is a ``bubble''
of concentration, by which we mean a region of spacetime of the form $\{ (t,x): |t-t_j| \leq c(\eta) N_j^{-2};
\quad |x-x_j| \leq c(\eta) N_j^{-1} \}$ inside the spacetime slab $I_j \times \R^n$ on which the solution $u$
has energy\footnote{Actually, we will only seek to obtain lower bounds on potential energy here, but
corresponding control on the kinetic energy can then be obtained by localized forms of the Sobolev inequality.}
 at least $c(\eta) > 0$.  Here $(t_j, x_j)$ is a point in $I_j \times \R^n$ and $N_j > 0$ is a frequency.
The spherical symmetry assumption allows us to choose $x_j = 0$; there is also a lower bound 
$N_j \geq c(\eta) |I_j|^{1/2}$ simply because the bubble has to be contained inside the slab $I_j \times \R^n$.
However, the harmonic analysis argument does not directly give an \emph{upper bound} on the frequency $N_j$; thus
the bubble may be much smaller than the slab.

In \cite{borg:scatter}, \cite{borg:book} an upper bound on $N_j$ is obtained by an \emph{induction on energy} 
argument; one assumes for contradiction that $N_j$ is very large, so the bubble is very small.  Without loss of 
generality we may assume the bubble lies in the lower half of the slab $I_j \times \R^n$.  Then when one evolves 
the bubble forward in time, it will have largely dispersed by the time it leaves $I_j \times \R^n$.  
Oversimplifying somewhat, the argument then proceeds by removing this bubble
(thus decreasing the energy by a non-trivial amount), applying an induction hypothesis to obtain Strichartz
bounds on the remainder of the solution, and then gluing the bubble back in by perturbation theory.  Unfortunately
it is this use of the induction hypothesis which eventually gives tower-exponential bounds rather than exponential
bounds in the final result.  Also there is some delicate playoff between various powers of $\eta$ which
needs additional care in four and higher dimensions.

Our main innovation is to obtain an upper bound on $N_j$ by more direct methods, dispensing with the need for
an induction on energy argument.  The idea is to use Duhamel's formula, to compare $u$ against the linear solutions
$u_\pm(t) := e^{i(t-t_\pm)\Delta} u(t_\pm)$.  We first eliminate a small number of intervals $I_j$ in which the linear
solutions $u_\pm$ have large $L^{2(n+2)/(n-2)}_{t,x}$ norm; the number of such intervals can be controlled by global
Strichartz estimates for the free (linear) Schr\"odinger equation.  Now let $I_j$ be one of the remaining intervals.  If
the bubble occurs in the lower half of $I_j$ then we\footnote{Again, this is an oversimplification; we must also dispose of
the non-linear interactions of $u$ with itself inside the interval $I_j$, but this can be done by some Strichartz
analysis and use of the pigeonhole principle.}  compare $u$ with $u_+$, taking advantage of the dispersive properties
of the propagator $e^{it\Delta}$ in our high-dimensional setting $n \geq 3$
to show that the error $u-u_+$ is in fact relatively smooth, which in turn implies
the bubble cannot be too small.  Similarly if the bubble occurs in the upper half of $I_j$ we compare $u$ 
instead with $u_-$.  Interestingly, there are some subtleties in very high dimension ($n \geq 6$) when the
non-linearity $F(u)$ grows quadratically or slower, as it now becomes rather difficult (in the large energy
setting) to pass from smallness of the non-linear solution (in spacetime norms)
to that of the linear solution or vice versa.

Once the bubble is shown to inhabit a sizeable portion of the slab, the rest of the argument essentially proceeds as
in \cite{borg:scatter}.  We wish to show that $J$ is bounded, so suppose for contradiction that $J$ is very large (so there are lots of bubbles).  
Then the Morawetz inequality \eqref{morawetz} can be used to show that 
the intervals $I_j$ must concentrate fairly rapidly at some point in time $t_*$; however one can then use localized mass
conservation laws to show that the bubbles inside $I_j$ must each shed a sizeable amount of mass (and energy) before
concentrating at $t_*$.  If $J$ is large enough there is so much mass and energy being shed that one can contradict
conservation of energy.  To put it another way, the mass conservation law implies that the bubbles cannot contract or expand rapidly, and the Morawetz inequality implies that the bubbles cannot persist stably for long periods of time.  Combining these two facts we can conclude that there are only a bounded number of bubbles.

It is worth mentioning that our argument is relatively elementary (compared against e.g. \cite{borg:scatter}, \cite{borg:book},
\cite{gopher}), especially in low dimensions $n=3,4,5$; the only tools are 
(non-endpoint) Strichartz estimates and Sobolev embedding, the Duhamel formula,
energy conservation, local mass conservation, and the Morawetz inequality, as well as some elementary combinatorial arguments. 
We do not need tools from Littlewood-Paley theory such as the para-differential 
calculus, although in the higher-dimensional cases $n \geq 6$ we will need fractional integration and the use of H\"older type estimates as a substitute for this para-differential calculus.  

\section{Notation and basic estimates}

We use $c, C > 0$ to denote various absolute constants depending only on the dimension $n$;
as we wish to track the dependence on the energy, we will \emph{not} allow these constants to depend on the
energy $E$.

For any time interval $I$, we use $L^q_t L^r_x(I \times \R^n)$ to denote the mixed spacetime Lebesgue norm
$$ \| u \|_{L^q_t L^r_x(I \times \R^n)} := (\int_I \| u(t) \|_{L^r(\R^n)}^q\ dt)^{1/q}$$
with the usual modifications when $q=\infty$.

We define the fractional differentiation operators $|\nabla|^\alpha := (-\Delta)^{\alpha/2}$ on $\R^n$.  Recall
that if $-n < \alpha < 0$ then these are fractional integration operators with an explicit form
\begin{equation}\label{fractint-def}
 |\nabla|^\alpha f(x) = c_{n,\alpha} \int_{\R^n} \frac{f(y)}{|x-y|^{n+\alpha}}\ dy
\end{equation}
for some computable constant $c_{n,\alpha} > 0$ whose exact value is unimportant to us; see e.g. \cite{stein:small}.
We recall that the Riesz transforms $\nabla |\nabla|^{-1} = |\nabla|^{-1} \nabla$ are bounded on $L^p(\R^n)$
for every $1 < p < \infty$; again see \cite{stein:small}.

\subsection{Duhamel's formula and Strichartz estimates}

Let $e^{it\Delta}$ be the propagator for the free Schr\"odinger equation $iu_t + \Delta u = 0$.  As is well known, this operator commutes with derivatives, and obeys the \emph{energy identity}
\begin{equation}\label{energy-identity}
\| e^{it\Delta} f \|_{L^2(\R^n)} = \| f\|_{L^2(\R^n)}
\end{equation}
and the \emph{dispersive inequality}
\begin{equation}\label{dispersive}
\| e^{it\Delta} f \|_{L^\infty(\R^n)} \leq C |t|^{-n/2} \| f\|_{L^1(\R^n)}
\end{equation}
for $t \neq 0$.  In particular we may interpolate to obtain the fixed-time estimates
\begin{equation}\label{dispersive-p}
\| e^{it\Delta} f \|_{L^{p}(\R^n)} \leq C |t|^{-n(\frac{1}{2} - \frac{1}{p})} \| f\|_{L^{p'}(\R^n)}
\end{equation}
for $2 \leq p \leq \infty$, where the dual exponent $p'$ is defined by $1/p + 1/p'$.

We observe \emph{Duhamel's formula}: if $iu_t + \Delta u = F$ on some time interval $I$,
then we have (in a distributional sense, at least)
\begin{equation}\label{duhamel}
 u(t) = e^{i(t-t_0)\Delta} u(t_0) - i \int_{t_0}^t e^{i(t-s)\Delta} F(s)\ ds
\end{equation}
for all $t_0, t \in I$, where we of course adopt the convention that $\int_{t_0}^t = - \int_t^{t_0}$
when $t < t_0$.  To estimate the terms on the right-hand side, we introduce the \emph{Strichartz norms}
$\dot S^k(I \times \R^n)$, defined for $k=0$ as
$$ \| u \|_{\dot S^0(I \times \R^n)} := \sup_{(q,r) \hbox{ admissible}} \| u \|_{L^q_t L^r_x(I \times \R^n)},$$
where admissibility was defined in \eqref{admissible}, and then for general\footnote{The homogeneous nature
of these norms causes some difficulties in interpreting elements of these spaces as a distribution
when $|k| \geq n/2$, but in practice we shall only work with $k=0,1$ and $n \geq 3$ and so these 
difficulties do not arise.} $k$ by
$$ \| u \|_{\dot S^k(I \times \R^n)} := \| |\nabla|^k u \|_{\dot S^0(I \times \R^n)}.$$
Observe that in the high dimensional setting $n \geq 3$, we have $2 \leq r < \infty$ for all
admissible $(q,r)$, so have boundedness of Riesz transforms (and thus we could replace $|\nabla|^k$ by
$\nabla^k$ for instance, when $k$ is a positive integer.  We note in particular that
\begin{equation}\label{strichartz-components}
\begin{split}
&\| \nabla^k u \|_{L^{2(n+2)/n}_{t,x}(I \times \R^n)}
+ \| \nabla^k u \|_{L^{2(n+2)/(n-2)}_t L^{2n(n+2)/(n^2 + 4)}_x(I \times \R^n)}\\&
+ \| \nabla^k u \|_{L^\infty_t L^2_x(I \times \R^n)} \leq C_k \|u\|_{\dot S^k(I \times \R^n)}
\end{split}
\end{equation}
for all positive integer $k \geq 1$.  Specializing further to the $k=1$ case we obtain
\begin{equation}\label{strichartz-components-2}
\| u \|_{L^{2(n+2)/(n-2)}_{t,x}(I \times \R^n)} +
\| u \|_{L^\infty_t L^{2n/(n-2)}_x(I \times \R^n)} \leq C \|u\|_{\dot S^1(I \times \R^n)}
\end{equation}
and in dimensions $n \geq 4$
\begin{equation}\label{strichartz-components-3}
\| u \|_{L^{2(n+2)/n}_t L^{2n(n+2)/(n^2 - 2n - 4)}_x(I \times \R^n)}
 \leq C \|u\|_{\dot S^1(I \times \R^n)}.
\end{equation}

We also define dual Strichartz spaces $\dot N^k(I \times \R^n)$, defined for $k=0$ as the Banach space dual of
$\dot S^0(I \times \R^n)$, and for general $k$ as
$$ \| F \|_{\dot N^k(I \times \R^n)} := \| |\nabla|^k F \|_{\dot N^0(I \times \R^n)}$$
(or equivalently, $\dot N^k$ is the dual of $\dot S^{-k}$).  From the first term in \eqref{strichartz-components}
and duality (and the boundedness of Riesz transforms) we observe in particular that
\begin{equation}\label{strichartz-dual}
\| F \|_{\dot N^k(I \times \R^n)} \leq \| \nabla^k F \|_{L^{2(n+2)/(n+4)}_{t,x}(I \times \R^n)}.
\end{equation}

We recall the \emph{Strichartz inequalities}
\begin{equation}\label{strichartz-homog}
 \| e^{i(t-t_0)\Delta} u(t_0) \|_{\dot S^k(I \times \R^n)} \leq C \| u(t_0) \|_{\dot H^k(\R^n)}
\end{equation}
and
\begin{equation}\label{strichartz-inhomog}
 \| \int_{t_0}^t e^{i(t-s)\Delta} F(s) \|_{\dot S^k(I \times \R^n)} \leq C \| F \|_{\dot N^k(I \times \R^n)};
\end{equation}
see e.g. \cite{tao:keel}; the dispersive inequality \eqref{dispersive-p} of course plays a key role in the
proof of these inequalities.  While we include the endpoint Strichartz pair $(q,r) = (2,\frac{2n}{n-2})$ in these estimates, this pair is not 
actually needed in our argument.  Observe that the constants $C$ here are independent of the choice of interval $I$.

\subsection{Local mass conservation}

We now recall a local mass conservation law appearing for instance in \cite{grillakis:scatter}; a related
result also appears in \cite{borg:scatter}.

Let $\chi$ be a bump function supported on the ball $B(0,1)$ which equals one on the ball $B(0,1/2)$
and is non-increasing in the radial direction.  For
any radius $R > 0$, we define the local mass $\Mass(u(t),B(x_0,R))$ of $u(t)$ on the ball $B(x_0,R)$ by
$$ \Mass(u(t),B(x_0,R)) := (\int \chi^2(\frac{x-x_0}{R}) |u(t,x)|^2\ dx)^{1/2};$$
note that this is a non-decreasing function of $R$.
Observe that if $u$ is a finite energy solution \eqref{nls}, then
$$ \partial_t |u(t,x)|^2 = - 2 \nabla_x \cdot \Im( \overline{u} \nabla_x u(t,x) )$$
(at least in a distributional sense), and so by integration by parts
$$ \partial_t \Mass(u(t),B(x_0,R))^2 = \frac{4}{R} \int \chi(\frac{x-x_0}{R}) (\nabla\chi)(\frac{x-x_0}{R}) 
\Im(\overline{u} \nabla_x u(t,x))\ dx$$
so by Cauchy-Schwarz
$$ |\partial_t \Mass(u(t),B(x_0,R))^2| \leq \frac{C}{R} \Mass(u(t),B(x_0,R)) (\int_{R/2 \leq |x-x_0| \leq R} |\nabla_x u(t,x)|^2\ dx)^{1/2}.$$
If $u$ has bounded energy $E(u) \leq E$, we thus have the approximate mass conservation law
\begin{equation}\label{mass-conserv}
 |\partial_t \Mass(u(t),B(x_0,R))| \leq C E^{1/2}/R.
\end{equation}
Observe that the same claim also holds if $u$ solves the free Schr\"odinger equation $iu_t + \Delta u = 0$
instead of the non-linear Schr\"odinger equation \eqref{nls}.  Note that the right-hand side decays with $R$.
This implies that if the local mass $\Mass(u(t),B(x_0,R))$ is large for some time $t$, then it can also
be shown to be similarly large for nearby times $t$, by increasing the radius $R$ if necessary to reduce
the rate of change of the mass.

From Sobolev and H\"older (or by Hardy's inequality) we can control the mass in terms of
the energy via the formula
\begin{equation}\label{hardy}
|\Mass(u(t), B(x_0,R))| \leq C E^{1/2} R.
\end{equation}

\subsection{Morawetz inequality}\label{review-sec}

We now give the proof of the Morawetz inequality \eqref{morawetz}; this inequality already
appears in \cite{borg:scatter}, \cite{borg:book}, \cite{grillakis:scatter} in three dimensions,
and the argument extends easily to higher dimensions, but for sake of completeness we give
the argument here.  

Using the scale invariance \eqref{scaling} we may rescale so that $A |I|^{1/2} = 1$.
We begin with the local momentum conservation identity
$$ 
\partial_t \Im(\partial_k u \overline{u})
= - 2\partial_j \Re(\partial_k u \overline{\partial_j u})
+ \frac{1}{2} \partial_k \Delta(|u|^2)
- \frac{2}{n-2} \partial_k |u|^{2n/(n-2)}$$
where $j, k$ range over spatial indices $1, \ldots, n$ with the usual summation conventions, and $\partial_k$
is differentiation with respect to the $x^k$ variable.  This identity can be verified directly from \eqref{nls};
observe that when $u$ is finite energy, both sides of this inequality make sense in the sense of distributions,
so this identity can be justified in the finite energy case by the local well-posedness theory\footnote{
For instance, one could smooth out the non-linearity $F$ (or add a parabolic dissipation term), obtain a similar law
for smooth solutions to the smoothed out equation, and then use the local well-posedness theory, see e.g. \cite{cwI}, 
to justify the process of taking limits.}.  If we multiply the above identity by the weight $\partial_k a$ for
some smooth, compactly supported  weight $a(x)$, and then integrate in space, we obtain (after some integration by parts)
\begin{align*}
 \partial_t \int_{\R^n} (\partial_k a) \Im(\partial_k u \overline{u})
=& 2\int_{\R^n} (\partial_j \partial_k a) \Re(\partial_k u \overline{\partial_j u})\\
&+ \frac{1}{2} \int_{\R^n} (-\Delta\Delta a) |u|^2\\
&+ \frac{2}{n-2} \int_{\R^n} \Delta a |u|^{2n/(n-2)}.
\end{align*}
We apply this in particular to the $C^\infty_0$ weight
$ a(x) := (\eps^2 + |x|^2)^{1/2} \chi(x)$,
where $\chi$ is a bump function supported on $B(0,2)$ which equals 1 on $B(0,1)$, and $0 < \eps < 1$ is a small
parameter which will eventually be sent to zero. In the region $|x| \leq 1$,
one can see from elementary geometry that $a$ is a convex function (its graph is a hyperboloid); in particular,
$(\partial_j \partial_k a) \Re(\partial_k u \overline{\partial_j u})$ is non-negative.  Further computation shows that
$$ \Delta a = \frac{n-1}{(\eps^2 + |x|^2)^{1/2}} + \frac{\eps^2}{(\eps^2 + |x|^2)^{3/2}}$$
and
$$ -\Delta \Delta a = \frac{(n-1)(n-3)}{(\eps^2 + |x|^2)^{3/2}} 
+ \frac{6(n-3) \eps^2}{(\eps^2 + |x|^2)^{5/2}}
+ \frac{15 \eps^4}{(\eps^2 + |x|^2)^{7/2}}$$
in this region; in particular $-\Delta \Delta a, \Delta a$ are positive in this region since $n \geq 3$.  
In the region $1 \leq |x| \leq 2$, $a$ and all of its derivatives are bounded uniformly in $\eps$, and
so the integrals here are bounded by $O(E(u))$ (using \eqref{hardy} to control the lower-order term).  
Combining these estimates we obtain the inequality
$$ \partial_t \int_{|x| \leq 2} (\partial_k a) \Im(\partial_k u \overline{u}) \geq
c \int_{|x| \leq 1} \frac{|u(t,x)|^{2n/(n-2)}}{(\eps^2 + |x|^2)^{1/2}}\ dx - C E(u).$$
Integrating this in time on $I$, and then using the fundamental theorem of calculus and the observation that
$a$ is Lipschitz, we obtain
$$ \sup_{t \in I} \int_{|x| \leq 2} |\nabla u(t,x)| |u(t,x)|\ dx \geq
c \int_I \int_{|x| \leq 1} \frac{|u(t,x)|^{2n/(n-2)}}{(\eps^2 + |x|^2)^{1/2}}\ dx - C E(u) |I|.$$
By \eqref{hardy} and Cauchy-Schwarz the left-hand side is $O(E(u))$.  Since $|I| = A^{-2} < 1$, we thus obtain
$$ \int_I \int_{|x| \leq 1} \frac{|u(t,x)|^{2n/(n-2)}}{(\eps^2 + |x|^2)^{1/2}}\ dx \leq C E(u).$$
Taking $\eps \to 0$ and using monotone convergence, \eqref{morawetz} follows.

\begin{remark}  In \cite{gopher}, an interaction variant of this Morawetz
inequality is used (superficially similar to the Glimm interaction potential 
as used in the theory of conservation laws), in which the weight $1/|x|$ is not present.  
In principle this allows for arguments such as the one here to
extend to the non-radial setting.  However the (frequency-localized)
interaction Morawetz inequality in \cite{gopher} is currently
restricted to three dimensions, and has a less favorable numerology\footnote{In the notation
of Corollary \ref{half-norm}, the interaction inequality in \cite{gopher} would give a bound of the form
$\sum_{I_j \subseteq I} |I_j|^{3/2} \leq C(\eta) (\max_{I_j \subseteq I} |I_j|)^{3/2}$,
which is substantially weaker and in particular does not seem to easily give the conclusions in Corollary \ref{burp}
or Proposition \ref{cc}, because the exponent $3/2$ here is greater than 1, whereas the corresponding
exponent $1/2$ arising from \eqref{morawetz} is less than 1.}
than \eqref{morawetz}, so it seems that the arguments given here
are insufficient to close the argument in the general case in higher dimensions.  
At the very least it seems that one would need to use more sophisticated control
on the movement of mass across frequency ranges, as is done in \cite{gopher}.
\end{remark}

\section{Proof of Theorem \ref{main}}

We now give the proof of Theorem \ref{main}.  The spherical symmetry of
$u$ is used in only one step, namely in Corollary \ref{concentration-cor},
to ensure that the solution concentrates at the spatial origin instead of at some other
location.

We fix $E$, $[t_-, t_+]$, $u$.  We may assume that the energy is large, $E > c > 0$, otherwise
the claim follows from the small energy theory.  From the bounded energy of $u$ we observe
the bounds
\begin{equation}\label{energy-bounded}
 \| u(t) \|_{\dot H^1(\R^n)} + \| u(t) \|_{L^{2n/(n-2)}(\R^n)} \leq C E^C
\end{equation}
for all $t \in [t_-, t_+]$.

We need some absolute constants $1 \ll C_0 \ll C_1 \ll C_2$, depending only on $n$, to be chosen later;
we will assume $C_0$ to be sufficiently large depending on $n$, $C_1$ sufficiently large depending on
$C_0, n$, and $C_2$ sufficiently large depending on $C_0, C_1, n$.
We then define the quantity $\eta := C_2^{-1} E^{-C_2}$.  Our task is to show that
$$ \int_{t_-}^{t_+} \int_{\R^n} |u(t,x)|^{2(n+2)/(n-2)}\ dx dt \leq C(C_0,C_1,C_2) \exp( C(C_0,C_1,C_2) E^{C(C_0,C_1,C_2)} ).$$
We may assume of course that
$$ \int_{t_-}^{t_+} \int_{\R^n} |u(t,x)|^{2(n+2)/(n-2)}\ dx dt > 4\eta$$
since our task is trivial otherwise.  We may then (by the greedy algorithm) subdivide $[t_-,t_+]$ 
into a finite number of disjoint intervals $I_1, \ldots, I_J$ for some
$J \geq 2$ such that
\begin{equation}\label{l10-small}
 \eta \leq \int_{I_j} \int_{\R^n} |u(t,x)|^{2(n+2)/(n-2)}\ dx dt \leq 2\eta
\end{equation}
for all $1 \leq j \leq J$.  It will then suffice to show that 
$$J \leq C(C_0,C_1,C_2) \exp( C(C_0,C_1,C_2) E^{C(C_0,C_1,C_2)} ).$$

We shall now prove various concentration properties of the solution on these intervals.  We begin
with a standard Strichartz estimate that bootstraps control on \eqref{l10-small} to control on all the
Strichartz norms (but we lose the gain in $\eta$):

\begin{lemma}\label{strichartz-control}  For each interval $I_j$ we have
$$ \| u \|_{\dot S^1(I_j \times \R^n)} \leq CE^C.$$
\end{lemma}

\begin{proof}  From Duhamel \eqref{duhamel}, Strichartz \eqref{strichartz-homog}, \eqref{strichartz-inhomog}
and the equation \eqref{nls} we have
$$ \| u \|_{\dot S^1(I_j \times \R^n)} \leq C \| u(t_j) \|_{\dot H^1(\R^n)} + \| F(u) \|_{\dot N^1(I_j \times \R^n)}$$
for any $t_j \in I_j$.  From \eqref{energy-bounded}, \eqref{strichartz-dual} we thus have
$$ \| u \|_{\dot S^1(I_j \times \R^n)} \leq C E^C + \| \nabla F(u) \|_{L^{2(n+2)/(n+4)}(I_j \times \R^n)}.$$
But from the chain rule and H\"older we have (formally, at least)
\begin{align*}
\| \nabla F(u) \|_{L^{2(n+2)/(n+4)}(I_j \times \R^n)} &\leq
C \| |u|^{4/(n-2)} |\nabla u| \|_{L^{2(n+2)/(n+4)}(I_j \times \R^n)} \\
&\leq C \| u \|_{L^{2(n+2)/(n-2)}_{t,x}(I_j \times \R^n)}^{4/(n-2)}
\|\nabla u \|_{L^{2(n+2)/n}(I_j \times \R^n)} \\
&\leq C \eta^{2/(n+2)} \| u \|_{\dot S^1(I_j \times \R^n)}
\end{align*}
by \eqref{l10-small}, \eqref{strichartz-components}.  Thus we have the formal inequality
$$  \| u \|_{\dot S^1(I_j \times \R^n)} \leq C E^C + C \eta^{2/(n+2)} \| u \|_{\dot S^1(I_j \times \R^n)}.$$
If $\eta$ is sufficiently small (by choosing $C_2$ large enough), then the claim follows, at least
formally.  To make the argument rigorous one can run a Picard iteration scheme that converges to the solution $u$
(see e.g. \cite{cwI} for details) and obtain the above types of bounds uniformly at all stages of the iteration; 
we omit the standard details.
\end{proof}

Next, we obtain \emph{lower} bounds on linear solution approximations to $u$ on an interval
where the $L^{2(n+2)/(n-2)}_{t,x}$ norm is small but bounded below.

\begin{lemma}\label{technical}  Let $[t_1, t_2] \subseteq [t_-, t_+]$ be an interval such that
\begin{equation}\label{u-dock}
 \eta/2 \leq \int_{t_1}^{t_2} \int_{\R^n} |u(t,x)|^{2(n+2)/(n-2)}\ dx dt \leq 2\eta.
\end{equation}
Then, if we define $u_l(t,x) := e^{i(t-t_l)\Delta}u(t_l)$ for $l=1,2$, we have
$$ \int_{t_1}^{t_2} \int_{\R^n} |u_l(t,x)|^{2(n+2)/(n-2)}\ dx dt \geq c \eta^C$$
for $l=1,2$.
\end{lemma}

\begin{proof}  Without loss of generality it suffices to prove the claim when $l=1$.  In low dimensions $n=3,4,5$
the Lemma is easy; indeed an inspection of the proof of Lemma \ref{strichartz-control} reveals that we have
the additional bound
$$ \| u - u_1 \|_{\dot S^1([t_1,t_2] \times \R^n)} \leq C E^C \eta^{2/(n+2)}$$
and hence by \eqref{strichartz-components-2}
$$ \| u - u_1 \|_{L^{2(n+2)/(n-2)}_{t,x}([t_1,t_2] \times \R^n)} \leq C E^C \eta^{2/(n+2)}.$$
When $n=3,4,5$ we have $2/(n+2) > (n-2)/2(n+2)$, and so the above estimates then show that
$u-u_1$ is smaller than $u$ in $L^{2(n+2)/(n-2)}_{t,x}([t_1,t_2] \times \R^n)$ norm
if $\eta$ is sufficienty small (i.e. $C_2$ is sufficiently large), at which point
the claim follows from the triangle inequality (and we can even replace $\eta^C$ by $\eta$).

In higher dimensions $n \geq 6$, the above simple argument breaks down. In fact
the argument becomes considerably more complicated (in particular, we were only able to
obtain a bound of $\eta^C$ rather than the more natural $\eta$); 
the difficulty is that while the non-linearity still decays faster than linearly as $u \to 0$, 
one of the factors is ``reserved'' for the derivative $\nabla u$, for which we have no smallness estimates, and
the remaining terms now decay linearly or worse, making it difficult to perform a perturbative analysis.
The resolution of this difficulty is rather technical, so we defer the proof of the higher dimensional
case to an Appendix (Section \ref{appendix}) so as not to interrupt the flow of the argument.
We remark however that the argument does not require any spherical symmetry assumption on the solution.
\end{proof}

Define the linear solutions $u_-$, $u_+$ on $[t_-,t_+] \times \R^n$ by
$u_\pm(t) := e^{i(t-t_\pm)\Delta} u(t_\pm)$; these are the analogue of
the scattering solutions for this compact interval $[t_-, t_+]$.
From \eqref{energy-bounded} and the Strichartz estimate \eqref{strichartz-homog}, \eqref{strichartz-components-2}, 
we have
$$ \int_{t_-}^{t_+} \int_{\R^n} |u_\pm(t,x)|^{2(n+2)/(n-2)}\ dx dt \leq C E^C.$$
Call an interval $I_j$ \emph{exceptional} if we have
$$ \int_{I_j} \int_{\R^n} |u_\pm(t,x)|^{2(n+2)/(n-2)}\ dx dt > \eta^{C_1}$$
for at least one choice of sign $\pm$, and \emph{unexceptional} otherwise.  
From the above global Strichartz estimate we see that
there are at most $O(E^C/\eta^{C_1})$ exceptional intervals, which will be acceptable for us from definition
of $\eta$.  Thus we may assume that there is at least one unexceptional interval.

Unexceptional intervals will be easier to control than exceptional ones, because the homogeneous component of Duhamel's formula
\eqref{duhamel} is negligible, leaving only the inhomogeneous component to be considered.  But as we shall see, this component
enjoys some additional regularity properties.  In particular, we now prove a concentration property of the solution on unexceptional intervals.

\begin{proposition}\label{concentration}  Let $I_j$ be an unexceptional interval.  Then there exists an $x_j\in \R^n$ such that
$$ \Mass( u(t), B(x_j, C \eta^{-C} |I_j|^{1/2}) ) \geq c \eta^{C C_0} |I_j|^{1/2}$$
for all $t \in I_j$.
\end{proposition}

\begin{proof} By time translation invariance and scale invariance \eqref{scaling} we may assume that
$I_j = [0,1]$.  We subdivide $I_j$ further into $[0,1/2]$ and $[1/2,1]$.  By \eqref{l10-small} and 
the pigeonhole principle
and time reflection symmetry if necessary we may assume that
\begin{equation}\label{upper-dens}
  \int_{1/2}^{1} \int_{\R^n} |u(t,x)|^{2(n+2)/(n-2)}\ dx dt > \eta/2.
\end{equation}
Since $I_j$ is unexceptional, we have
\begin{equation}\label{unexceptional}
 \int_0^1 \int_{\R^n} |u_-(t,x)|^{2(n+2)/(n-2)}\ dx dt \leq \eta^{C_1}.
\end{equation}
By \eqref{unexceptional}, \eqref{l10-small} and the pigeonhole principle, we may find an interval 
$[t_* - \eta^{C_0}, t_*] \subset [0,1/2]$ such that\footnote{In the low dimensional
case $n=3,4,5$ we may skip this pigeonhole step.  Indeed from \eqref{upper-dens}, \eqref{unexceptional}
 and Duhamel we may conclude that $\int_{t_-}^{0} e^{i(t-s)\Delta} F(u(s))\ ds$ has
large $L^{2(n+2)/(n-2)}_{t,x}$ norm on the slab $[1/2,1] \times \R^n$; this is because the proof of 
Lemma \ref{technical} shows that the effect of the forcing terms arising from the time interval $[0,1]$
are of size $O(\eta^{4/(n-2)})$, which is smaller than $\eta/2$ for $n = 3,4,5$; one then continues the proof from
\eqref{largeness-2} onwards with only minor changes.  However this simple
argument does not seem to work in higher dimensions.}
\begin{equation}\label{pigeonhole}
  \int_{t_* - \eta^{C_0}}^{t_*} \int_{\R^n} |u(t,x)|^{2(n+2)/(n-2)}\ dx dt < C \eta^{C_0}.
\end{equation}
and
\begin{equation}\label{unexceptional-frozen}
 \int_{\R^n} |u_-(t_* - \eta^{C_0},x)|^{2(n+2)/(n-2)}\ dx \leq C \eta^{C_1}.
\end{equation}
Applying Lemma \ref{technical} to the time interval
$[t_*, 1]$ we see that
\begin{equation}\label{largeness}
 \int_{t_*}^1 \int_{\R^n} |(e^{i(t-t_*)\Delta}u(t_*))(x)|^{2(n+2)/(n-2)}\ dx dt \geq c \eta^C.
\end{equation}
By Duhamel's formula \eqref{duhamel} we have
\begin{equation}\label{duh}
e^{i(t-t_*)\Delta} u(t_*) = u_-(t) 
- i \int_{t_* - \eta^{C_0}}^{t_*} e^{i(t-s)\Delta} F(u(s))\ ds
- i \int_{t_-}^{t_* - \eta^{C_0}} e^{i(t-s)\Delta} F(u(s))\ ds.
\end{equation}
Since $I_j$ is unexceptional, we have
$$ \int_{t_*}^1 \int_{\R^n} |u_-(t,x)|^{2(n+2)/(n-2)}\ dx dt \leq \eta^{C_1}.$$
From \eqref{pigeonhole} and Lemma \ref{strichartz-control}, it is easy to see (using the chain rule and
H\"older as in the proof of Lemma \ref{strichartz-control}) that
\begin{equation}\label{fun}
  \| F(u) \|_{\dot N^1( [t_* - \eta^{C_0}, t_*] \times \R^n )} \leq C E^C \eta^{cC_0},
\end{equation}
and hence by Strichartz \eqref{strichartz-inhomog}
$$ \int_{t_*}^1 \int_{\R^n} |\int_{t_* - \eta^{C_0}}^{t_*} e^{i(t-s)\Delta} F(u(s))\ ds|^{2(n+2)/(n-2)}(x)\ dx dt
\leq C E^C \eta^{cC_0}.$$
From these estimates and \eqref{largeness}, we thus see from the triangle inequality
(if $C_0$ is large enough, and $\eta$ small enough (i.e. $C_2$ large enough depending on $C_0$)) that
\begin{equation}\label{largeness-2}
\| v \|_{L^{2(n+2)/(n-2)}_{t,x}([t_*,1] \times \R^n)}
 \geq c \eta^C
\end{equation}
where $v$ is the function
\begin{equation}\label{v-def}
v := \int_{t_-}^{t_* - \eta^{C_0}} e^{i(t-s)\Delta} F(u(s))\ ds.
\end{equation}
We now complement this lower bound on $v$ with an upper bound.  First observe from Lemma \ref{strichartz-control}
that 
$$ \| u \|_{\dot S^1([t_*,1] \times \R^n)} \leq C E^C;$$
also from \eqref{energy-bounded} and \eqref{strichartz-homog} we have
$$ \| u_- \|_{\dot S^1([t_*,1] \times \R^n)} \leq C E^C.$$
Finally, from \eqref{fun} and \eqref{strichartz-inhomog}
$$ \| \int_{t_* - \eta^{C_0}}^{t_*} e^{i(t-s)\Delta} F(u(s))\ ds \|_{\dot S^1([t_*,1] \times \R^n)} \leq C E^C.$$
From the triangle inequality and \eqref{duh} we thus have
\begin{equation}\label{v-s}
\| v \|_{\dot S^1([t_*,1] \times \R^n)} \leq C E^C.
\end{equation}
We shall need some additional regularity control on $v$.
For any $h \in \R^n$, let $u^{(h)}$ denote the translate of $u$ by $h$, i.e. $u^{(h)}(t,x) := u(t,x-h)$.

\begin{lemma}\label{hold}  We have the bound
$$ \| v^{(h)} - v \|_{L^\infty_t L^{2(n+2)/(n-2)}_x([t_*,1] \times \R^n)} \leq C E^C \eta^{-CC_0} |h|^c$$
for all $h \in \R^n$. 
\end{lemma}

\begin{proof}  First consider the high-dimensional case $n \geq 4$.  We use \eqref{energy-bounded},
the chain rule and H\"older to observe that
\begin{align*}
 \| \nabla F(u(s)) \|_{L^{2n/(n+4)}(\R^n)} &\leq C
\| |u(s)|^{4/(n-2)} |\nabla u(s)| \|_{L^{2n/(n+4)}(\R^n)}\\
&\leq C \| u(s) \|_{L^{2n/(n-2)}(\R^n)}^{4/(n-2)} \| \nabla u(s) \|_{L^2(\R^n)} \\
&\leq C E^C,
\end{align*}
so by the dispersive inequality \eqref{dispersive-p} 
$$ \| \nabla e^{i(t-s)\Delta} F(u(s)) \|_{L^{2n/(n-4)}(\R^n)} \leq C E^C |t-s|^{-2}.$$
Integrating this for $s$ in $[t_-,t_* - \eta^{C_0}]$ we obtain 
$$ \| \nabla v \|_{L^\infty_t L^{2n/(n-4)}_x([t_*,t_1] \times \R^n)} \leq C E^C \eta^{-CC_0};$$
interpolating this with \eqref{v-s}, \eqref{strichartz-components} we obtain
$$ \| \nabla v \|_{L^\infty_t L^{2(n+2)/(n-2)}_x([t_*,t_1] \times \R^n)} \leq C E^C \eta^{-CC_0}.$$
The claim then follows (with $c=1$) from the Fundamental theorem of calculus and Minkowski's inequality.

Now consider the three-dimensional case $n=3$.  From \eqref{energy-bounded}, the fundamental theorem
of calculus, and Minkowski's inequality we have
$$ \| u^{(h)}(s) - u(s) \|_{L^2(\R^3)} \leq C E^C |h|,$$
while from the triangle inequality we have
$$ \| u^{(h)}(s) - u(s) \|_{L^6(\R^3)} \leq C E^C,$$
and hence
$$ \| u^{(h)}(s) - u(s) \|_{L^3(\R^3)} \leq C E^C |h|^{1/2}.$$
Since $F(u)$ is quintic in three dimensions, we thus have from H\"older and \eqref{energy-bounded} that
\begin{align*}
\| F(u)^{(h)}(s) - F(u)(s) \|_{L^1(\R^3)} &\leq C \| |u^{(h)}(s) - u(s)| (|u^{(h)}(s)| + |u(s)|)^4 \|_{L^1(\R^3)}\\
&\leq C \| u^{(h)}(s) - u(s) \|_{L^3(\R^3)} \| u(s) \|_{L^6(\R^3)}^4 \leq C E^C |h|^{1/2}.
\end{align*}
Integrating this for $s \in [t_-, t_* - \eta^{C_0}]$ using \eqref{dispersive} we obtain 
$$ \| v^{(h)} - v \|_{L^\infty_{t,x}([t_*,1] \times \R^n)} \leq C E^C \eta^{-CC_0} |h|^{1/2}.$$
On the other hand, from \eqref{v-s}, \eqref{strichartz-components-2}, and the triangle inequality we have
$$ \| v^{(h)} - v \|_{L^\infty_t L^6_x([t_*,1] \times \R^n)} \leq C E^C \eta^{-CC_0}$$
and the claim follows by interpolation.
\end{proof}

We can average this lemma over all $|h| \leq r$, for some scale $0 < r < 1$ to be chosen shortly, to obtain
$$ \| v_{av} - v \|_{L^\infty_t L^{2(n+2)/(n-2)}_x([t_*,1] \times \R^n)} \leq C E^C \eta^{-CC_0} r^c$$
where $v_{av}(x) := \int \chi(y) v(x + ry)\ dy$ for some bump function $\chi$ supported on $B(0,1)$ 
of total mass one.  In particular by
a H\"older in time we have
$$ \| v_{av} - v \|_{L^{2(n+2)/(n-2)}_{t,x}([t_*,1] \times \R^n)} \leq C E^C \eta^{-CC_0} r^c.$$
Thus if we choose $r := \eta^{CC_0}$ for some large enough $C$, and $\eta$ is sufficiently small, we see
from \eqref{largeness-2} that
$$ \| v_{av} \|_{L^{2(n+2)/(n-2)}_{t,x}([t_*,1] \times \R^n)} \geq c \eta^C.$$
On the other hand, by H\"older and Young's inequality
\begin{align*}
 \| v_{av} \|_{L^{2n/(n-2)}_{t,x}([t_*,1] \times \R^n)}
&\leq C \| v_{av} \|_{L^\infty_t L^{2n/(n-2)}_x([t_*,1] \times \R^n)}\\
&\leq C \| v \|_{L^\infty_t L^{2n/(n-2)}_x([t_*,1] \times \R^n)} \\
&\leq C E^C
\end{align*}
by \eqref{v-s}, \eqref{strichartz-components}.  Thus by H\"older we have
$$ \| v_{av} \|_{L^\infty_{t,x}([t_*,1] \times \R^n)} \geq c \eta^C E^{-C}.$$
Thus we may find a point $(t_j, x_j) \in [t_*,1] \times \R^n$ such that
$$ |\int \chi(y) v(t_j, x_j+ry)\ dy| \geq c \eta^C E^{-C},$$
and in particular by Cauchy-Schwarz
$$ \Mass(v(t_j), B(x_j,R)) \geq c \eta^C E^{-C} r^C.$$
for all $R \geq r$.  Observe from \eqref{v-def}
that $v$ solves the free Schr\"odinger equation on $[t_* - \eta^{C_0}, 1]$,
and has energy $O(E^C)$ by \eqref{v-s}, \eqref{strichartz-components}.  Thus by
\eqref{mass-conserv} we have
$$ \Mass(v(t_* - \eta^{C_0}), B(x_j,R)) \geq c \eta^C E^{-C} r^C$$
for all $t \in [t_*, 1]$, if we set $R := C \eta^{-C} E^C r^{-C}$ for some appropriate constants $C$.
From Duhamel's formula \eqref{duhamel} (or \eqref{duh}) we have 
$$ u(t_* - \eta^{C_0}) = u_-(t_* - \eta^{C_0}) - iv(t_* - \eta^{C_0}).$$
From \eqref{unexceptional-frozen} and H\"older we have
$$ \Mass(u_-(t_* - \eta^{C_0}), B(x_j,R)) \leq C R^C \eta^{CC_1}.$$
Thus if we choose $C_1$ sufficiently large depending on $C_0$ (recalling that $r = \eta^{CC_0}$ and
$R = C \eta^{-C} E^C r^{-C}$), and assume $\eta$ sufficiently small depending polynomially on $E$, we have
$$ \Mass(u(t_* - \eta^{C_0}), B(x_j,R)) \geq c \eta^C E^{-C} r^C.$$
By another application of \eqref{mass-conserv} we thus have
$$ \Mass(u(t), B(x_j,\eta^{-CC_0})) \geq c \eta^{-CC_0}$$
for all $t \in [0,1]$, and Proposition \ref{concentration} follows.
\end{proof}

We now exploit the radial symmetry of $u$ to place the concentration point $x_j$
at the origin.  This is the only place where the spherical symmetry assumption is used.

\begin{corollary}\label{concentration-cor}  Let $I_j$ be an unexceptional interval, and assume that the solution $u$ is spherically symmetric.  Then 
we have
$$ \Mass( u(t), B(0, C \eta^{-CC_0} |I_j|^{1/2}) ) \geq c \eta^{CC_0} |I_j|^{1/2}$$
for all $t \in I_j$.
\end{corollary}

\begin{proof}  We again rescale $I_j = [0,1]$.
Let $x_j$ be as in Proposition \ref{concentration}.  Fix $t \in [0,1]$.
If $|x_j| = O(\eta^{-C' C_0})$ for some
$C'$ depending only on $n$ then we are done.  Now suppose that $|x_j| \geq \eta^{-C' C_0}$.  Then if
$C'$ is big enough, we can find $\eta^{-cC'}$ rotations of the ball $B(x_j, C \eta^{-C C_0})$
which are disjoint.  On each one of these balls, the mass of $u(t)$ is at least $c \eta^{C C_0}$
by the spherical symmetry assumption; by H\"older this shows that the $L^{2n/(n-2)}$ norm
of $u(t)$ on these balls is also $c \eta^{C C_0}$.  Adding this up for each of the $\eta^{-cC' C_0}$ balls,
we obtain a contradiction to \eqref{energy-bounded} if $C' C_0$ is large enough.  Thus
we have $|x_j| = O(\eta^{-C' C_0})$ and the claim follows.
\end{proof}

From this corollary and H\"older we see that
$$ \int_{|x| \leq R} \frac{|u(t,x)|^{2n/(n-2)}}{|x|}\ dx dt \geq c
\eta^{C C_0} |I_j|^{-1/2}$$
whenever $t \in I_j$ for some unexceptional interval $I_j$,
and $R \geq C \eta^{-C C_0} |I_j|^{1/2}$.  In particular
we have
$$ \int_{I_j} \int_{|x| \leq R} \frac{|u(t,x)|^{2n/(n-2)}}{|x|}\ dx dt \geq c
\eta^{C C_0} |I_j|^{1/2}.$$

Combining this with \eqref{morawetz} and the bounded energy we obtain the following
combinatorial bound on the distribution of the intervals $I_j$.

\begin{corollary}\label{half-norm}  Assume that the solution $u$ is spherically symmetric.
For any interval $I \subseteq [t_-, t_+]$, we have
$$ \sum_{1 \leq j \leq J: I_j \subseteq I} |I_j|^{1/2} \leq C \eta^{-C(C_0,C_1)} |I|^{1/2}.$$
(note we can use $\eta^{-C}$ to absorb any powers of the energy which appear; also, note
that the $O(C \eta^{-C_1})$ exceptional intervals cause no difficulty).
\end{corollary}

This bound gives quite strong control on the possible distribution of the intervals $I_j$, for instance we have

\begin{corollary}\label{burp} Assume that the solution $u$ is spherically symmetric.
  Let $I = \bigcup_{j_1 \leq j \leq j_2} I_j$ be a union of consecutive intervals.  Then there exists
$j_1 \leq j \leq j_2$ such that $|I_j| \geq c \eta^{C(C_0,C_1)} |I|$.
\end{corollary}

\begin{proof}  From the preceding corollary we have
$$ C \eta^{-C(C_0,C_1)} |I|^{1/2} \geq \sum_{j_1 \leq j \leq j_2} |I_j|^{1/2}
\geq \sum_{j_1 \leq j \leq j_2} |I_j| (\sup_{j_1 \leq j \leq j_2} |I_j|)^{-1/2}.$$
Since $\sum_{j_1 \leq j \leq j_2} |I_j| = |I|$, the claim follows.
\end{proof}

We now repeat a combinatorial argument\footnote{It seems of interest to remove the logarithm in
this Proposition, since this would make our final estimate polynomial in the energy instead
of exponential.  It seems however one cannot achieve this purely on the strength of the
Morawetz estimate \eqref{morawetz} and the mass conservation law \eqref{mass-conserv}, as the control
on the intervals $I_j$ provided by these two estimates does not preclude the possibility for the energy to concentrate 
on a Cantor set of times of dimension less than $1/2$, which can use up an exponential number of intervals before
the local mass conservation begins to conflict with energy conservation.  One possibility is to combine the Morawetz
inequality \eqref{morawetz} with the interaction Morawetz inequalities in \cite{gopher}, although those inequalities
are in some sense even weaker and thus less able to control the total number of intervals. We remark that for the cubic NLS in three
dimensions, the known bounds are polynomial in the energy and mass \cite{ckstt:french}, \cite{ckstt:cubic-scatter}, but this
is because the equation is $H^1$-subcritical and $L^2$-supercritical, which force the lengths $|I_j|$ of the intervals
to be bounded both above and below.  See \cite{nak:scatter} for a related discussion.}
of Bourgain \cite{borg:scatter} to show that the intervals $I_j$
must now concentrate at some time $t_*$:

\begin{proposition}\label{cc}  Assume that the solution $u$ is spherically symmetric.
Then there exists a time $t_* \in [t_-, t_+]$ and distinct unexceptional intervals $I_{j_1}, \ldots, I_{j_K}$
for some $K > c \eta^{C(C_0,C_1)} \log J$ such that
\begin{equation}\label{geom-decay}
 |I_{j_1}| \geq 2 |I_{j_2}| \geq 4 |I_{j_3}| \geq \ldots \geq 2^{K-1} |I_{j_K}|
\end{equation}
and such that
$ \dist(t_*, I_{j_k}) \leq C \eta^{-C(C_0,C_1)} |I_{j_k}|$
for all $1 \leq k \leq K$.
\end{proposition}

\begin{proof}  We run the algorithm from Bourgain \cite{borg:scatter}.  We first recursively define a nested sequence of
intervals $I^{(k)}$, each of which is a union of consecutive unexceptional $I_j$, as follows.  We first remove the
$O(\eta^{-C_1})$ exceptional intervals from $[t_-, t_+]$, leaving $O(\eta^{-C_1})$ connected components.
One of these, call it $I^{(1)}$, must be the union of $J_1 \geq c \eta^{C_1} J$ consecutive unexceptional 
intervals.  By Corollary \ref{burp}, there
exists an $I_{j_1} \subseteq I^{(1)}$ such that $|I_{j_1}| \geq c \eta^{C C_0} |I^{(1)}|$, so in particular
$\dist(t, I_{j_1}) \leq C \eta^{-C C_0} |I_{j_1}|$ for all $t \in |I^{(1)}|$.  Now we remove $I_{j_1}$ from $I^{(1)}$,
and more generally remove all intervals $I_j$ from $I^{(1)}$ for which $|I_j| > |I_{j_1}|/2$.  There can
be at most $C \eta^{-C C_0}$ such intervals to remove, since $I_{j_1}$ was so large.  If $J_1 \leq C \eta^{-C C_0}$
then we set $K = 1$ and terminate the algorithm.  Otherwise, we observe that the remaining connected components
of $I^{(1)}$ still contain at least $c \eta^{C C_0} J$ intervals, and there are $O(\eta^{-C C_0})$ such components.  Thus by
the pigeonhole principle we can find one of these components, $I^{(2)}$, which is the union of $J_2 \geq c \eta^{C C_0} J_1$ intervals, 
each of which must have length less than or equal to $|I_{j_1}|/2$ by construction.
Now we iterate the algorithm, using Corollary \ref{burp} to locate an interval $I_{j_2}$ in $I^{(2)}$
such that $|I_{j_2}| \geq c \eta^{C C_0} |I^{(2)}|$, and then removing all intervals of length $> |I_{j_2}|/2$
from $I^{(2)}|$.  If the number of intervals in $|I^{(2)}|$ is $O(\eta^{-C C_0})$, we terminate the algorithm,
otherwise we can pass as before to a smaller interval $I^{(3)}$ which is a union of $J_3 \geq c\eta^{C C_0} J_2$
intervals.  We can continue in this manner for $K$ steps for some $K > c \eta^{C(C_0,C_1)} \log J$ until we run out
of intervals.  The claim then follows by choosing $t_*$ to be an arbitrary time in $I^{(K)}$.
\end{proof}

Let $t_*$ and $I_{j_1}, \ldots, I_{j_k}$ be as in the above Proposition.
From Proposition \ref{concentration} we recall that
$$ \Mass(u(t), B(x_{j_k}, C \eta^{-C(C_0)} |I_{j_k}|^{1/2})) \geq c \eta^{C(C_0)} |I_{j_k}|^{1/2}$$
for all $t \in I_{j_k}$.  Applying \eqref{mass-conserv} and adjusting the constants $c, C$ as necessary we thus see
that
$$ \Mass(u(t_*), B_k) \geq c \eta^{C(C_0,C_1)} |I_{j_k}|^{1/2},$$
where each $B_k$ is a ball $B_k := B(x_{j_k}, C \eta^{-C(C_0,C_1)} |I_{j_k}|^{1/2})$.  On the other hand, from \eqref{hardy}
we observe that
$$ \Mass(u(t_*), B_k) \leq C \eta^{-C(C_0,C_1)} |I_{j_k}|^{1/2}.$$
Let $N := C_2 \log (1/\eta)$.  If we choose this constant $C_2$ large enough, we thus see from the above mass bounds
and \eqref{geom-decay} that
$$ \sum_{k+N \leq k' \leq K} \int_{B_{k'}} |u(t_*,x)|^2\ dx \leq \frac{1}{2}
\int_{B_k} |u(t_*,x)|^2\ dx,$$
and hence
$$ \int_{B_k \backslash (\bigcup_{k+N \leq k' \leq K} B_{k'})} |u(t_*,x)|^2 \geq c \eta^{C(C_0,C_1)} |I_{j_k}|.$$
Applying H\"older's inequality\footnote{An alternate approach here is to use the spherical symmetry to
move the balls to be centered at the origin, and apply Hardy's inequality, see \cite{borg:scatter}, \cite{grillakis:scatter}.  
However this approach shows that one does not need the spherical symmetry assumption to conclude the argument
provided that one has a concentration result similar to Proposition \ref{cc}.}, we thus obtain
$$ \int_{B_k \backslash (\bigcup_{k+N \leq k' \leq K} B_{k'})} |u(t_*,x)|^{2n/(n-2)} \geq c \eta^{C(C_0,C_1)}.$$
Summing this in $k$ and telescoping, we obtain
$$ \int_{\R^n} |u(t_*,x)|^{2n/(n-2)} \geq c \eta^{C(C_0,C_1)} K / N.$$
Using \eqref{energy-bounded} we thus obtain
$$ K \leq C \eta^{-C(C_0,C_1)} N E^{C} \leq C(C_0,C_1,C_2) \eta^{-C(C_0,C_1)}.$$
Since $K > c \eta^{C(C_0,C_1)} \log J$, we obtain 
$$J \leq \exp(C \eta^{-C(C_0,C_1)}) \leq \exp(C(C_0,C_1,C_2) E^{C(C_0,C_1,C_2)})$$
as desired. This proves Theorem \ref{main}.

\begin{remark}\label{remark-gopher}  One can use Proposition \ref{concentration} to improve the bounds obtained in
\cite{gopher} in the non-radial case, as one no longer needs to use the induction hypothesis
to obtain concentration bounds on the solution.  It may also be possible to use a variant of the
techniques here to also obtain the reverse Sobolev inequality.  However the remaining portion of
the arguments seem to require a heavier use of the induction hypothesis (in order to obtain certain
frequency localization properties of the energy), and so we were unable to fully remove the tower-type bounds
from the result in \cite{gopher}.
\end{remark}

\section{Appendix: Proof of Lemma \ref{technical} in high dimensions}\label{appendix}

We now give the rather technical proof of Lemma \ref{technical} in the high-dimensional
case $n \geq 6$; the idea is to find an iteration scheme which converges acceptably after
the first few terms, leaving us to estimate a finite number of iterates (which we can estimate by
more inefficient means).
We differentiate \eqref{nls} and use the chain rule to obtain the equation
$$ (i \partial_t + \Delta) \nabla u = V_1 \nabla u + V_2 \overline{\nabla u}$$
where $V_1 := \frac{n}{n-2} |u|^{\frac{4}{n-2}}$ and $V_2 := \frac{2}{n-2}
|u|^{\frac{4}{n-2}} \frac{u^2}{|u|^2}$.  From \eqref{u-dock} we have
\begin{equation}\label{V-bound}
 \| V_1 \|_{L^{(n+2)/2}_{t,x}([t_1,t_2] \times \R^n)} + \| V_2 \|_{L^{(n+2)/2}_{t,x}([t_1,t_2] \times \R^n)} 
\leq C \eta^c,
\end{equation}
which by \eqref{strichartz-dual}, \eqref{strichartz-components}, and H\"older implies in particular
that
\begin{equation}\label{contraction}
 \| V_1 w + V_2 \overline{w} \|_{\dot N^0([t_1,t_2] \times \R^n)} \leq
C \eta^c \| w \|_{\dot S^0([t_1,t_2] \times \R^n)}.
\end{equation}
From Duhamel's formula \eqref{duhamel} we have $\nabla u = \nabla u_1 + A \nabla u$, where $A$ is the (real) linear operator
$$ Aw(t) := -i \int_{t_1}^t e^{i(t-s)\Delta}(V_1(s) w(s) + V_2(s) \overline{w}(s))\ ds.$$
From Strichartz \eqref{strichartz-inhomog} and \eqref{contraction} we see that 
\begin{equation}\label{A-bound}
\| A w \|_{\dot S^0([t_1,t_2] \times \R^n)}
\leq C \eta^c \| w \|_{\dot S^0([t_1,t_2] \times \R^n)};
\end{equation}
thus for $\eta$ sufficiently small, $A$ is a contraction on $\dot S^0([t_1,t_2] \times \R^n)$.  
Also, from Strichartz \eqref{strichartz-homog} and \eqref{energy-bounded}
we see that 
\begin{equation}\label{u1-bound}
\| \nabla u_1 \|_{\dot S^0([t_1,t_2] \times \R^n)} \leq C E^C.
\end{equation}
  Thus for some absolute constant $M$ (depending only on $n$), we see that
we have the Neumann series approximation
$$ \| \nabla u - \sum_{m=0}^M A^m \nabla u_1 \|_{\dot S^0([t_1,t_2] \times \R^n)} \leq \eta$$
(for instance), assuming that $\eta$ is sufficiently small depending (polynomially) on the energy. 
Now introduce the spacetime norm
$$ \| w \|_X := \| |\nabla|^{-1} w \|_{L^{2(n+2)/(n-2)}_{t,x}([t_1,t_2] \times \R^n)}$$
where $|\nabla|^{-1} := (-\Delta)^{-1/2}$.  From \eqref{strichartz-components-2} (and the boundedness
of Riesz transforms) we observe that
\begin{equation}\label{fractint}
 \| w \|_X \leq C \| w \|_{\dot S^0([t_1,t_2] \times \R^n)}
\end{equation}
and hence
$$ \| \nabla u - \sum_{m=0}^M A^m \nabla u_1 \|_X \leq C\eta.$$
On the other hand, from \eqref{u-dock} and Calder\'on-Zygmund theory we have
$$ \| \nabla u \|_X \geq c \eta^{(n-2)/2(n+2)}.$$
Thus by the triangle inequality, we have
$$ \| A^m \nabla u_1 \|_X \geq c \eta^{(n-2)/2(n+2)}$$
for some $0 \leq m \leq M$, again assuming that $\eta$ is sufficiently small.

Ideally we would now like the operator $A$ to be bounded on $X$.  We do not know
if this is true; however we have the following weaker (and technical) version of this fact which
suffices for our application.

\begin{lemma}\label{tech}  For any $w \in \dot S^0([t_1,t_2] \times \R^n)$, we have the estimate 
$\| A w \|_X \leq C E^C \| w \|_{\dot S^0([t_1,t_2] \times \R^n)}^{1-\theta} \| w \|_X^\theta$
for some absolute constant $0 < \theta < 1$ (depending only on $n$).
\end{lemma}

Assuming this Lemma for the moment, we apply it together with \eqref{A-bound}, \eqref{u1-bound}, 
\eqref{fractint} we obtain a bound of the form
$$ \| A^m \nabla u_1 \|_X \leq C_m \eta^{-C_m} E^{C_m} \| \nabla u_1 \|_X^{\theta_m}$$
for some constants $C_m, \theta_m > 0$.  Combining this with our lower bound on $\|A^m \nabla u_1 \|_X$ we
obtain $\| \nabla u_1 \|_X \geq c \eta^C$ (assuming $\eta$ sufficiently small depending on $E$, and allowing
constants to depend on the fixed constant $M$), and Lemma \ref{technical} follows (again using the boundedness of
Riesz transforms).

It remains to prove Lemma \ref{tech}.  The point is to take advantage of one of the (many) refinements 
of the
Sobolev embedding used\footnote{More precisely, we need a statement to the effect that
the Sobolev theorem is only sharp if one of the ``wavelet coefficients'' of the function is extremely
large (close to its maximal size).  The argument below could be reformulated as an interpolation
inequality (of Gagliardo-Nirenberg type) for Triebel-Lizorkin spaces, but we have elected to give a direct
argument that does not rely on too much external machinery.}  to prove \eqref{fractint}; we shall
use an argument based on Hedberg's inequality.  We will not attempt to gain powers of $\eta$ here 
(since the Neumann series step has in some sense fully exploited those gains already) and
so shall simply discard all such gains that we encounter.

We will make the \emph{a priori} assumption that $w$ is smooth and rapidly decreasing;
this can be removed by the usual limiting argument.
We normalize $\| w\|_{\dot S^0([t_1,t_2] \times \R^n)} := 1$, and write $\alpha := \| w \|_X$,
thus $\alpha \leq C$ by \eqref{fractint}.  Our task is to show that $\| Aw \|_X \leq C \alpha^c$.

Observe from \eqref{A-bound} and \eqref{strichartz-components-3} that
$$ \| |\nabla|^{-1} A w \|_{L^{2(n+2)/n}_t L^{2n(n+2)/(n^2 - 2n - 4)}_x([t_1,t_2] \times \R^n)}
\leq C$$
and hence by H\"older it will suffice to show that
$$ \| |\nabla|^{-1} A w(t) \|_{L^{2n/(n-2)}(\R^n)}
\leq C \alpha^c$$
for all $t \in [t_1,t_2]$.

By time translation invariance we may set $t=0$.  Write $v := Aw(0)$.  
We now use a variant of Hedberg's inequality.  From \eqref{A-bound} and \eqref{strichartz-components} we have
\begin{equation}\label{l2-bounded}
 \| v \|_{L^2_x} \leq C;
\end{equation}
if we let $M$ denote the Hardy-Littlewood maximal function
$$ Mv(x) := \sup_{r > 0} \frac{1}{|B(x,r)|} \int_{B(x,r)} |v(y)|\ dy,$$
then by the Hardy-Littlewood maximal inequality (see \cite{stein:small}) we have
$$ \| (M v)^{n-2/n} \|_{L^{2n/(n-2)}_x}
= \| Mv \|_{L^2_x}^{(n-2)/n} \leq C.$$
It thus suffices to prove the pointwise Hedberg-type inequality
$$ ||\nabla|^{-1} v|(x) \leq C \alpha^c (Mv(x))^{(n-2)/n}.$$
We may translate so that $x = 0$.  Set $R := (Mv(0))^{-2/n}$ (the case $Mv(0) = 0$ being trivial),
thus
$$ \int_{B(0,r)} |v(y)|\ dy \leq C R^{-n/2} r^n$$
for all $r > 0$.  From \eqref{l2-bounded} and Cauchy-Schwarz we also have
$$ \int_{B(0,r)} |v(y)|\ dy \leq C r^{n/2}$$
and thus
\begin{equation}\label{v-ball}
 \int_{B(0,r)} |v(y)|\ dy \leq C r^{n/2} (1 + \frac{r}{R})^{-n/2}.
\end{equation}
By \eqref{fractint-def}, it suffices to show that
$$ |\int_{\R^n} \frac{v(y)}{|y|^{n-1}}\ dy| \leq C \alpha^c R^{1 - \frac{n}{2}}.$$
By the above estimates, we see that we can prove this estimate in the regions
$|y| \leq \alpha^{c_0} R$ and $|y| \geq \alpha^{-c_0} R$, even if we place the absolute values
inside the integral, where $0 < c_0\ll 1$ is an absolute
constant to be chosen shortly.  Thus it will suffice to estimate the remaining region
$\alpha^{c_0} R \leq |y| \leq \alpha^{-c_0} R$.  Partitioning the integral via smooth cutoffs, we see
that it suffices (if $c_0$ was chosen sufficiently small) to show that
$$ |\int_{\R^n} v(y) \varphi(y/r)\ dy| \leq C \alpha^c r^{\frac{n}{2}},$$
for all $r > 0$, where $\varphi$ is a real-valued bump function.  One may verify from dimensional analysis
that this estimate (as well as the hypotheses) are invariant under the scaling
$$ w(t,x) \mapsto \lambda^{-n/2} w(t/\lambda^2, x/\lambda); \quad
v(t,x) \mapsto \lambda^{-n/2} w(t/\lambda^2, x/\lambda)$$
and so we may take $r = 1$.  Since $v = Aw(0)$, we thus reduce to proving that
$$ |\langle Aw(0), \varphi \rangle| \leq C \alpha^c.$$
Expanding out the definition of $A$ and using duality, we can write this as
\begin{equation}\label{duel}
 |\int_{t_1}^0 \int_{\R^n} (V_1(t) w(t) + V_2(t) \overline{w}(t)) e^{it\Delta} \varphi\ dx dt| \leq C \alpha^c.
\end{equation}
From \eqref{V-bound}, \eqref{strichartz-components} we have
$$ \| V_1 w + V_2 \overline{w} \|_{L^{2(n+2)/(n+4)}_{t,x}([t_1,0] \times \R^n)} \leq C,$$
while a direct computation\footnote{Indeed, one just needs to note that $e^{-it\Delta} \varphi$ is bounded in $L^2_x$ and decays
in $L^\infty_x$ like $O(t^{-n/2})$ to verify this claim.} of $e^{-it\Delta} \varphi$ shows that
$$ \| e^{-it\Delta} \varphi \|_{L^{2(n+2)/n}_{t,x}([t_1,-\tau] \times \R^n)} \leq C \tau^{-c}$$
for all $\tau > 1$.  Thus if we set $\tau = C \alpha^{-c_0}$ for some small $c_0$ to be chosen later,
we see that the portion of \eqref{duel} arising from $[t_1,-\tau]$ is acceptable, and it suffices to then
prove the bound on $[-\tau,0]$.  In fact we will prove the fixed time estimates
$$
 |\int_{\R^n} (V_1(t) w(t) + V_2(t) \overline{w(t)}) e^{it\Delta} \varphi \ dx| \leq C E^C (1+t)^{C}
\| |\nabla|^{-1} w(t) \|_{L^{2(n+2)/(n-2)}(\R^n)}^c
$$
for all $t \in [-\tau,0]$, which proves the claim if $c_0$ is sufficiently small, thanks to H\"older's inequality
and the hypothesis $\|w\|_X = \alpha$.

Fix $t$.  We shall just prove this inequality for $V_2 \overline{w}$, as the corresponding estimate for $V_1 w$ is 
similar.  Because of the negative derivative on $w$ on the right-hand side, we shall need some regularity control 
on $V_2$.  Note that $V_2$ behaves like $|u|^{4/(n-2)}$; since $4/(n-2) \leq 1$,
the standard fractional chain rule is not easy to apply.  Instead, we will work in H\"older-type spaces\footnote{Using 
H\"older spaces rather than Sobolev spaces costs an epsilon of regularity (see e.g. \cite{stein:small} for a discussion)
but for our purposes any non-zero amount of regularity will suffice.  The reader may recognize the arguments below as that of splitting a product into paraproducts; however we are avoiding the use of standard 
paraproduct theory as it does not interact well with non-linear maps such as $u \mapsto V_2$ which may only be H\"older continuous of order
$4/(n-2) < 1$.},
which are more elementary.  As with Lemma \ref{hold}, we let $u^{(h)}$ denote the translate $u^{(h)}(x) := u(x-h)$ of $u$ by $h$ for any $h \in \R^n$;
similarly define $V_2^{(h)}$, etc.  From \eqref{energy-bounded}, the fundamental theorem of calculus, and Minkowski's
inequality we have
$$ \| u^{(h)}(t) - u(t) \|_{L^2(\R^n)} \leq C E^C |h|.$$
Since the function $z \mapsto |z|^{\frac{4}{n-2}} \frac{z^2}{|z|^2}$ is H\"older continuous of order $4/n-2$,
we thus have the pointwise inequality
$$ |V_2^{(h)}(t) - V_2(t)| \leq C |u^{(h)} - u|^{4/(n-2)} $$ 
which gives the H\"older type bounds
$$ \| V_2^{(h)}(t) - V_2(t) \|_{L^{(n-2)/2}(\R^n)} \leq C E^C |h|^{4/(n-2)}.$$
From \eqref{strichartz-components} and the normalization $\|w\|_{\dot S^0([t_1,t_2] \times \R^n)} = 1$ we have
$ \| w(t) \|_{L^2_x} \leq C$,
and hence by H\"older (and the decay of $e^{-it\Delta} \varphi$ in space)
$$ |\int_{\R^n} (V_2(t) - V_2^{(h)}) \overline{w(t)} e^{it\Delta} \varphi \ dx| \leq C E^C (1+t)^{C} |h|^{4/(n-2)}.$$
Similarly, from \eqref{energy-bounded} we have
\begin{equation}\label{v2-bound}
 \| V_2^{(h)}(t) \|_{L^{n/2}(\R^n)} \leq C
\| u(t) \|_{L^{2n/(n-2)}(\R^n)}^{4/(n-2)} \leq C E^C;
\end{equation}
a direct computation also shows that
$$ \| e^{it\Delta}\varphi - (e^{it\Delta}\varphi)^{(h)} \|_{L^{2n/(n-4)}(\R^n)} \leq C (1+t)^C |h| \leq C (1+t)^C |h|^{4/(n-2)}$$
for $|h| \leq 1$ (say), and so by H\"older again
$$ |\int_{\R^n} V_2^{(h)} \overline{w(t)} (e^{it\Delta} \varphi - (e^{it\Delta}\varphi)^{(h)})\ dx| 
\leq C E^C (1+t)^{C} |h|^{4/(n-2)}.$$
Combining this with the previous estimate we obtain
$$ |\int_{\R^n} \overline{w(t)} (V_2 e^{it\Delta} \varphi - 
(V_2 e^{it\Delta} \varphi)^{(h)})\ dx| \leq C E^C (1+t)^{C} |h|^{4/(n-2)},$$
or equivalently that
$$ |\int_{\R^n} (\overline{w(t)} - \overline{w(t)}^{(-h)}) V_2 e^{it\Delta} \varphi\ dx| \leq 
C E^C (1+t)^{C} |h|^{4/(n-2)}.$$
We can average this over all $|h| \leq r$, where the radius $0 < r < 1$ will be chosen later, to obtain
\begin{equation}\label{w-diff}
 |\int_{\R^n} (\overline{w(t)} - \overline{w_{av}(t)}) V_2 e^{it\Delta} \varphi\ dx| \leq 
C E^C (1+t)^{C} r^{4/(n-2)}
\end{equation}
where $w_{av}(t,x) := \int \chi(y) w(t,x+r y)$ for some bump function $\chi$ of total mass 1.
On the other hand, from the H\"ormander multiplier theorem (see \cite{stein:small}) and some Fourier
analysis we see that
$$ \| w_{av}(t) \|_{L^{2(n+2)/(n-2)}(\R^n)} \leq C r^{-C} \| |\nabla|^{-1} w(t) \|_{L^{2(n+2)/(n-2)}(\R^n)},$$
and by combining this with \eqref{v2-bound} and decay estimates on $e^{it\Delta} \varphi$ we obtain
$$ |\int_{\R^n} \overline{w_{av}(t)} V_2 e^{it\Delta} \varphi\ dx| \leq 
C E^C (1+t)^{C} r^{-C} \| |\nabla|^{-1} w(t) \|_{L^{2(n+2)/(n-2)}(\R^n)}.$$
Combining this with \eqref{w-diff} and optimizing in $r$ we obtain Lemma \ref{tech}, and thus Lemma \ref{technical}.

\end{document}